\documentclass[12pt]{article}

\usepackage {amssymb}
\usepackage {amscd}
\usepackage {latexsym}
\usepackage {euscript}

\usepackage[T2A]{fontenc}
\usepackage[cp1251]{inputenc}
\usepackage[english,russian]{babel}
\usepackage[tbtags]{amsmath}
\usepackage{amsfonts,amssymb}

\usepackage{mathrsfs}
\usepackage{graphicx}

\newtheorem{lemma}{Lemma}
\newtheorem{theorem}{Theorem}
\newtheorem{definition}{Definition}
\newtheorem{proposition}{Proposition}

\newtheorem{remark}{Remark}

\begin{document}
{\Large
\begin{center} Sharp embeddings of uniformly localized\\ Bessel potential spaces into multiplier spaces
 \end{center}}

\bigskip
\centerline{A.~A.Belyaev and A.~A.~Shkalikov}\footnote{This work is supported by Russian Science Foundation (RNF) under grant No 17-11-01215.}

\bigskip
\medskip

\centerline{{\bf Abstract}}
 \medskip
 For $p > 1, \: \gamma \in \mathbb R$, denote by $H^{\gamma}_p(\mathbb{R}^n)$ the Bessel potential space, by $H^{\gamma}_{p, \; unif}(\mathbb{R}^n)$ the corresponding uniformly localized Bessel potential space and by $M[s,-t]$ the space of multipliers from $H^s_2(\mathbb{R}^n)$ into  $H^{-t}_2(\mathbb{R}^n)$. Assume that $s, \: t \geqslant 0, \: n/2 > \max(s, \: t) > 0, \: r: = \min(s, \: t), \: p_1: = n/max(s, \: t)$. Then following embeddings hold
$$
H^{-r}_{p_1, \: unif}(\mathbb{R}^n) \subset M[s, -t] \subset H^{-r}_{2, \: unif}(\mathbb{R}^n).
$$
The main result of the paper claims the sharpness of the left embedding
in the following sense: it does not hold if the lower index $p_1$ is replaced by $p_1 -\varepsilon$ with any small $\varepsilon > 0$.

\medskip

{\bf Key words:} Bessel  potential spaces, multipliers, embedding theorems.

\bigskip

{ \Large  1. Introduction}


\medskip

In this paper we are concerned with constructive description of multiplier space from one Bessel potential space $H^s_2(\mathbb{R}^n)$ into another such space $H^{-t}_2(\mathbb{R}^n)$ for $s, \; t \geqslant 0$. Namely, our main goal is to establish the sharp character of some parameters naturally arising in the problem of finding such a constructive description in terms of uniformly localized Bessel potential spaces.

Describing multiplier space in these terms is of great importance for the singular perturbation theory of many particular differential operators. In a model case of Laplace operator, defined on the classical Sobolev space $W^2_2(\mathbb{R}^n)$, study of its perturbations by singular potentials is closely related to the description problem for multipliers acting from the Sobolev space with positive smoothness index into the Sobolev space with negative smoothness index. Indeed, classical Laplace operator can be extended to an operator acting from $W^1_2(\mathbb{R}^n)$ to $W^{-1}_2(\mathbb{R}^n)$ and this operator is correctly defined if and only if perturbation potential is a multiplier from $W^1_2(\mathbb{R}^n)$ to $W^{-1}_2(\mathbb{R}^n)$. Moreover, employing multiplier technique for more general Bessel potential spaces turned out to be extremely useful not only in the case of Laplace operator's singular perturbations but also for general strongly elliptic operators (see \cite{NZSh2}).

One of the main tools in the study of perturbed operator's spectral properties is the use of resolvent convergence in order to approximate initial perturbed operator of general type by operators perturbed by specific potentials. This approach has a long, fruitful history as both more general approach with a setting in abstract Hilbert space (see, e.g., early works on topic \cite{Stum, Greenlee} and also \cite{VBorisov90} for applications of this abstract approach to some specific perturbations) and numerous applications to concrete differential operators were developed. In the latter direction various techniques were used to obtain results of this type: for example, in \cite{Daners03} resolvent convergence was employed to study convergence of Dirichlet problem solutions on sequence of domains to the solution of this problem on a limit domain, while other contributions were motivated by problems in mathematical physics, see, for instance, \cite{BrFiTe}, where weak convergence of Radon measures was employed in order to approximate Schr\"odinger operator $-\Delta + \mu$ for some finite Radon measure $\mu$ by Hamiltonians describing point interactions. The most important motivation for us lies in the fact that under natural assumptions norm convergence of perturbations in the multiplier space yields uniform resolvent convergence of corresponding perturbations. In \cite[Theorem 7]{NZSh2} there was established a general result for the strongly elliptic operator, perturbed by non--discrete singular distribution from Bessel potential space with negative smoothness index, when uniform resolvent approximation by perturbations with smooth potentials is considered.

Therefore, finding a description of multiplier space in terms of its coincidence with some specific space from the scale $H^s_{p, \: unif}(\mathbb{R}^n)$, introduced implicitly by R.\,S.~Strichartz in \cite{Str67}, gives us constructive conditions on perturbation potential in order to develop spectral theory for perturbed operators even in a case when singular perturbation potential neither has a discrete support nor is a specific potential (e.g., Coulomb potential or some other extensively treated cases, for references see \cite{AGRH2ndEd}). For Bessel potential spaces $s_0 = \frac{n}{p}$ serves as a limit case for the problem of finding a description of the type
$$
M[H^s_p(\mathbb{R}^n) \to H^s_p(\mathbb{R}^n)] = H^s_{p, \: unif}(\mathbb{R}^n), \: s \geqslant 0, \: p > 1,
$$
obtained in \cite{Str67}. Moreover, the situation is analogous for more general spaces of Besov and Lizorkin--Triebel type as their descriptions in terms of uniformly localized Besov--Lizorkin--Triebel spaces can be obtained under similar assumptions (see \cite{Fr, Sickel99} and monograph \cite[Chapter IV]{RSbook} for references). When we consider the multiplier space $M[H^s_p(\mathbb{R}^n) \to H^t_p(\mathbb{R}^n)]$ with positive smoothness indices $s$ and $t$ no longer being equal but satisfying natural assumption $s \geqslant t$ (otherwise multiplier space is trivial), the role of limiting values still belongs to $s_0 = \frac{n}{p}$.

Yet this is no longer a case for the multipliers acting between Bessel potential spaces with smoothness indices of different sign. This situation was treated in a series of papers by A.\,A.~Shkalikov, M.~I. Neiman--Zade, J.\,G.~Bak and A.\,A.~Belyaev (see \cite{NZSh1, BSh, NZSh2, BelShk, Bel2}). A slight reformulation of \cite[Lemma 4, Lemma 6]{NZSh2} states that if $s, \: t \geqslant 0$ and $s + t \neq 0$, then
$$
M[H^s_2(\mathbb{R}^n) \to H^{-t}_2(\mathbb{R}^n)] = H^{-\min(s, \; t)}_{2, \: unif}(\mathbb{R}^n),
$$
whenever $\max(s, \; t) > \frac{n}{2}$, while continuous embeddings
\begin{equation}\label{main_embeddings}
H^{-\min(s, \: t)}_{p_1, \: unif}(\mathbb{R}^n) \subset M[H^s_2(\mathbb{R}^n) \to H^{-t}_2(\mathbb{R}^n)] \subset H^{-\min(s, \; t)}_{2, \: unif}(\mathbb{R}^n),
\end{equation}
where $p_1 = \frac{n}{\max(s, \; t)}$, are valid whenever $\max (s, t) < \frac{n}{2}.$

In this paper we prove that a parameter $p_1 = \frac{n}{\max(s, \; t)}$ in the embedding \eqref{main_embeddings} is sharp in a following sense: for arbitrary  $\varepsilon$, satisfying conditions
$$
0 < \varepsilon < \frac{n}{\max (s, t)} - 2,
$$
there exists a distribution
$$
u_{\varepsilon} \in H^{- \min (s, t)}_{p_1 - \varepsilon, \: unif}(\mathbb{R}^n) \setminus \: M[H^s_2(\mathbb{R}^n), \: H^{-t}_2(\mathbb{R}^n)].
$$
The main importance of this result lies in the fact that it guarantees absence of the constructive description for $M[H^s_2(\mathbb{R}^n) \to H^{-t}_2(\mathbb{R}^n)]$ in terms of the uniformly localized Bessel potential spaces.

\bigskip
\bigskip
\bigskip
\bigskip
\bigskip

{\Large 2. Basic definitions and classical results.}

\bigskip

\bigskip

In what follows, we say that a Banach space $X$ is continuously embedded in a Banach space $Y$ and write $X \subset Y$ if $X$ is embedded in $Y$ in a set--theoretic sense and the estimate
$$
\|u\|_Y \leqslant C \cdot \|u\|_X \; \; \forall \: u \in X
$$
is valid for a constant $C > 0$ independent of $u$.

By $D(\mathbb{R}^n)$ and $D'(\mathbb{R}^n)$ we denote the space of smooth compactly supported test functions and the corresponding dual space of distributions (generalized functions), while by $S(\mathbb{R}^n)$ and $S'(\mathbb{R}^n)$ we mean Schwartz space of rapidly decreasing functions and dual Schwartz space of tempered distributions respectively.

For arbitrary $p \in (1; \; +\infty)$ we denote by $p'$ its Lebesgue dual, i.e. such a number $p'$ that
$$
\frac{1}{p} + \frac{1}{p'} = 1.
$$

\begin{definition}\label{def_H^0_p}
For $p \geqslant 1$, we say that the distribution $\mathbf{v} \in D'(\mathbb{R}^n)$ belongs to the space $H^0_p(\mathbb{R}^n)$ if there exists $v \in L_p(\mathbb{R}^n)$ such that
$$
\mathbf{v}(f) = \int\limits_{\mathbb{R}^n} f \cdot v \; d\mu \; \; \forall \: f \in D(\mathbb{R}^n),
$$
and the norm of $\mathbf{v}$ in the space $H^0_p(\mathbb{R}^n)$ is defined by
$$
\|\mathbf{v}\|_{H^0_p(\mathbb{R}^n)} = \|v\|_{L_p(\mathbb{R}n)}.
$$
\end{definition}

\begin{definition}\label{def_H^s_p}
Let $s \in \mathbb{R}, \; p \geqslant 1$. Let an operator $J_s \colon S'(\mathbb{R}^n) \to S'(\mathbb{R}^n)$ be defined by
$$
J_s(u) \stackrel{def}{=} \mathcal{F}^{-1}(\varphi_s \cdot \mathcal{F}(u)) \; \; \forall \: u \in S'(\mathbb{R}^n),
$$
where $\mathcal{F}$ and $\mathcal{F}^{-1}$ are the direct and inverse Fourier transforms in the dual Schwartz space $S'(\mathbb{R}^n)$ and $\varphi_s(x) = (1 + |x|^2)^{\frac{s}{2}}$ for all $x \in \mathbb{R}^n$. Then Bessel potential space $H^s_p(\mathbb{R}^n)$ is defined as a set of all distributions $u \in S'(\mathbb{R}^n)$ such that $J_s(u) \in H^0_p(\mathbb{R}^n)$, equipped with the norm
$$
\|u\|_{H^s_p(\mathbb{R}^n)} = \|J_s(u)\|_{H^0_p(\mathbb{R}^n)}.
$$
\end{definition}

We note that for arbitrary $s \in \mathbb{R}$ and $p > 1$ the space $(H^s_p(\mathbb{R}^n))^{*}$, defined as a dual space to $H^s_p(\mathbb{R}^n)$, is isometrically isomorphic to the space $H^{-s}_{p'}(\mathbb{R}^n)$.

\begin{remark}\label{multiplication_operator}
Let us consider $C^{\infty}_b(\mathbb{R}^n)$, i.e. the space of infinitely differentiable functions bounded together with all their derivatives on $\mathbb{R}^n$, and assume that $s \in \mathbb{R}, \: p \geqslant 1$ and $f \in C^{\infty}_b(\mathbb{R}^n)$. Then the operator $A_f \colon H^s_p(\mathbb{R}^n) \to H^s_p(\mathbb{R}^n)$ of multiplication by a function $f$ is well--defined and bounded with respect to the norm $\| \cdot \|_{H^s_p(\mathbb{R}^n)}$. In particular, the operator of multiplication by a function $\varphi \in D(\mathbb{R}^n)$ is well--defined as a continuous operator on $H^s_p(\mathbb{R}^n)$ and the norm of this operator does not change if the function $\varphi$ is replaced by a function $\varphi_{(z)}, \: z \in \mathbb{R}^n$, defined as
$$
\varphi_{(z)}(x) = \varphi(x - z) \; \; \forall \: x \in \mathbb{R}^n.
$$
\end{remark}

In what follows we use the following generalization of classical Sobolev embedding theorem for the Bessel potential spaces.

{\bfseries Embedding theorem for Bessel potential spaces} {\itshape Let $p, \: q \geqslant 1$, and let $s, \: t \in \mathbb{R}$. If $p \leqslant q$ and $s - \frac{n}{p} \geqslant t - \frac{n}{q}$, then
$$
H^s_p(\mathbb{R}^n) \subset H^t_q(\mathbb{R}^n).
$$}

\begin{definition}\label{loc_unif_spaces}
Let $s \in \mathbb{R}, \; p \geqslant 1$. Then
$$
H^s_{p, \: loc}(\mathbb{R}^n) \stackrel{def}{=} \{u \in D'(\mathbb{R}^n) \: | \; f \cdot u \in H^s_p(\mathbb{R}^n) \; \; \forall \: f \in D(\mathbb{R}^n)\}
$$
and uniformly localized Bessel potential space is defined as follows:
$$
H^s_{p, \: unif, \: \eta}(\mathbb{R}^n) \stackrel{def}{=} \{ u \in H^s_{p, \: loc}(\mathbb{R}^n) \: | \; \|u\|_{s, \: p, \: unif, \: \eta} \stackrel{def}{=} \sup\limits_{z \in \mathbb{R}^n} \|\eta_{(z)} \cdot u\|_{H^s_p(\mathbb{R}^n)} < +\infty\},
$$
where $\eta_{(z)}(x) = \eta(x - z) \; \; \forall \; x \in \mathbb{R}^n$ and function $\eta \in D(\mathbb{R}^n)$ meets the conditions
$$
0 \leqslant \eta(x) \leqslant 1 \; \; \forall \: x \in \mathbb{R}^n, \; \;  \eta(x) = 1 \; \; \forall \: x \in \mathbb{R}^n \colon |x| \leqslant 1, \; \; \eta(x) = 0 \; \; \forall \: x \in \mathbb{R}^n \colon |x| \geqslant 2.
$$
\end{definition}

For different $\eta \in D(\mathbb{R}^n)$, satisfying the conditions of Definition \ref{loc_unif_spaces}, norms $\| \cdot \|_{s, \: p, \: unif, \eta}$ are equivalent to each other. Therefore, we omit the index $\eta$ and simply write $H^s_{p, \: unif}(\mathbb{R}^n)$ for uniformly localized Bessel potential space and $\| \cdot \|_{H^s_{p, \: unif}(\mathbb{R}^n)}$ for its norm.

With Remark \ref{multiplication_operator} taken into account, it is also easy to see that we have continuous embedding
$$
H^s_p(\mathbb{R}^n) \subset H^s_{p, \: unif}(\mathbb{R}^n).
$$

\begin{remark}\label{unif_embedding} Generalization of Sobolev embedding theorem yields continuous embedding
$$
H^s_{p, \: unif}(\mathbb{R}^n) \subset H^t_{q, \: unif}(\mathbb{R}^n),
$$
whenever conditions of the Sobolev embedding theorem are met. Moreover, it is easy to demonstrate that this continuous embedding also holds for $p\geqslant q \geqslant 1$ and $s \geqslant t, \: s, \: t \in \mathbb{R}$.
\end{remark}

\begin{definition}\label{def_multipliers}
Let $s, \: t \geqslant 0, \;  p, \: q > 1$. A distribution $\mu \in H^{-t}_{q', \: loc}(\mathbb{R}^n)$ belonging to the set
$$
\{\mu \in H^{-t}_{q', \: loc}(\mathbb{R}^n)\: | \; \exists \: C > 0 \colon  \|f \cdot \mu\|_{H^{-t}_{q'}(\mathbb{R}^n)} \leqslant C \cdot \|\mathbf{f}\|_{H^s_p(\mathbb{R}^n)} \; \; \forall \: f \in  D(\mathbb{R}^n)\}
$$
is called a multiplier from the space $H^s_p(\mathbb{R}^n)$ to the space $H^{-t}_{q'}(\mathbb{R}^n)$ and the norm on the multiplier space $M[H^s_p(\mathbb{R}^n) \to H^{-t}_{q'}(\mathbb{R}^n)]$ is introduced by
$$
\|\mu\|_{M[H^s_p(\mathbb{R}^n) \to H^{-t}_{q'}(\mathbb{R}^n)]} = \inf\{C > 0 \: | \; \|f \cdot \mu\|_{H^{-t}_{q'}(\mathbb{R}^n)} \leqslant C \|\mathbf{f}\|_{H^s_p(\mathbb{R}^n)} \; \; \forall \: f \in D(\mathbb{R}^n)\}.
$$
\end{definition}

Since the space $D(\mathbb{R}^n)$ is densely embedded in $H^s_p(\mathbb{R}^n)$ for $s \in \mathbb{R}, \: p \geqslant 1$, it follows that under assumptions of Definition \ref{def_multipliers}, a multiplier $\mu \in M[H^s_p(\mathbb{R}^n) \to H^{-t}_{q'}(\mathbb{R}^n)]$ uniquely defines a bounded linear operator $M_{\mu} \colon H^s_p(\mathbb{R}^n) \to H^{-t}_{q'}(\mathbb{R}^n)$, such that
$$
M_{\mu}(\mathbf{f}) = f \cdot \mu \; \; \forall \: f \in D(\mathbb{R}^n).
$$
Moreover, if $f \in D(\mathbb{R}^n)$ and $\mathbf{f} \in M[H^s_p(\mathbb{R}^n) \to H^{-t}_{q'}(\mathbb{R}^n)]$, then
$$
M_{\mathbf{f}}(u) = f \cdot u \; \; \forall \: u \in H^s_p(\mathbb{R}^n).
$$

With the remark on the isomorphism between $(H^t_q(\mathbb{R}^n))^{*}$ and $H^{-t}_{q'}(\mathbb{R}^n)$ taken into account, we obtain equivalent definition of the multiplier space $M[H^s_p(\mathbb{R}^n) \to H^{-t}_{q'}(\mathbb{R}^n)]$.

\begin{definition}\label{equiv_def_multi}
Let $s, \: t \geqslant 0, \; p, \: q > 1$. Distribution $\mu \in D'(\mathbb{R}^n)$ is a multiplier from $H^s_p(\mathbb{R}^n)$ to $H^{-t}_{q'}(\mathbb{R}^n)$, if there exists a constant $C > 0$, such that
$$
|\mu(f \cdot g)| \leqslant C \cdot \|\mathbf{f}\|_{H^s_p(\mathbb{R}^n)} \cdot \|\mathbf{g}\|_{H^t_q(\mathbb{R}^n)} \; \; \forall \: f, \: g \in D(\mathbb{R}^n),
$$
with the multiplier norm being defined by
$$
\| \mu\|_{M[H^s_p(\mathbb{R}^n) \to H^{-t}_{q'}(\mathbb{R}^n)]} =  \inf\{C > 0 \: | \; |\mu(f \cdot g)| \leqslant C \cdot \| \mathbf{f}\|_{H^s_p(\mathbb{R}^n)} \cdot \| \mathbf{g}\|_{H^t_q(\mathbb{R}^n)} \; \; \forall \: f, \: g \in D(\mathbb{R}^n)\}
$$
for all $\mu \in M[H^s_p(\mathbb{R}^n) \to H^{-t}_{q'}(\mathbb{R}^n)]$.
\end{definition}


This allows us to obtain the following technical result.

\begin{lemma}\label{multipliers_symmetry}
Let $s, \: t \geqslant 0$ and $p, \: q > 1$. Then
$$
M[H^s_p(\mathbb{R}^n) \to H^{-t}_{q'}(\mathbb{R}^n)] = M[H^t_q(\mathbb{R}^n) \to H^{-s}_{p'}(\mathbb{R}^n)]
$$
and norms of these spaces are equal.
\end{lemma}

This assertion immediately follows from Definition \ref{equiv_def_multi}.

The following result which served as a principal motivation for our research was obtained by A.\,A.~Shkalikov and M.\,I.~Neiman--Zade in 2006.

{\bfseries Theorem} (see \cite[Lemma 4, Lemma 5]{NZSh2}){\itshape Let $s, \: t \geqslant 0$ and $\max(s, \: t) \in \left(0; \: \frac{n}{2}\right)$. Then
$$
H^{-\min(s, \: t)}_{p_1, \: unif}(\mathbb{R}^n) \subset M[H^s_2(\mathbb{R}^n) \to H^{-t}_2(\mathbb{R}^n)] \subset H^{-\min(s, \: t)}_{2, \: unif}(\mathbb{R}^n),
$$
where $p_1 = \frac{n}{\max(s, \; t)}$.}

For the sake of brevity in what follows we shall use notation
$$
M[s, \: -t] \stackrel{def}{=} M[H^s_2(\mathbb{R}^n) \to H^{-t}_2(\mathbb{R}^n)].
$$

\bigskip

\bigskip

\bigskip

{\Large 3. A specific class of regular distributions: principal properties and belonging to uniformly localized Bessel potential space}

\bigskip

\medskip

In this section we consider the regular functional $\mathbf{f_{\alpha}}, \; \alpha > 0$, generated by a real--valued function
$$
f_{\alpha}(x) \stackrel{def}{=} \begin{cases}
|x|^{-\alpha}, \; \; x \in \mathbb{R}^n \setminus \{0\},\\
0, \; \; x = 0.
\end{cases}
$$
We examine conditions on $\alpha$ which guarantee that this functional is correctly defined as a tempered distribution, i.e. an element from the dual Schwartz space $S'(\mathbb{R}^n)$, and also find a necessary condition for $\mathbf{f_{\alpha}}$ to belong to the space $H^{-t}_{2, \: unif}(\mathbb{R}^n)$.

Converting to generalized spherical coordinates, it can be readily deduced that for arbitrary $r_1, \: r_2 > 0$
$$
\int\limits_{B(r_1, r_2)} f_{\alpha}(x) dx = C(n) \cdot \int\limits_{[r_1, r_2]} r^{-\alpha} \cdot r^{n - 1} dr = C(n) \cdot \int\limits_{[r_1, r_2]} r^{-(\alpha - n + 1)} dr,
$$
where constant
$$
C(n) = \frac{2 \cdot \pi^{\frac{n}{2}}}{\Gamma(\frac{n}{2})}
$$
is a the surface area of the $(n - 1)$--dimensional hypersphere. Hence, $f_{\alpha}$ is integrable on $B_1$ if and only if $\alpha < n$ and $f_{\alpha}$ is integrable on $\mathbb{R}^n \setminus B_1$ if and only if $\alpha > n$.

Furthermore, a simple technical check gives the following criterion for functional $\mathbf{f_{\alpha}}$ to belong to the dual Schwartz space $S'(\mathbb{R}^n)$.

\begin{proposition}\label{correctness_prop}
Let $\alpha > 0.$ Then a functional $\mathbf{f_{\alpha}} \colon S(\mathbb{R}^n) \to \mathbb{R}$, defined by
$$
\varphi \stackrel{\mathbf{f_{\alpha}}}{\longmapsto} \int\limits_{\mathbb{R}^n} f_{\alpha}(x) \cdot \varphi(x) \: d x,
$$
is well--defined and continuous on Schwartz space $S'(\mathbb{R}^n)$ if and only if $\alpha < n$.
\end{proposition}

Now we remind some technical propositions needed for the proof of this section's main result.

Firstly, we recall a well--known result from Riesz potential theory.

\begin{proposition}\label{Fourier_transform}
$(see, e.g., $\: \cite[\S5.1, Lemma 1]{Stein}$)$ Let $0 < \alpha < n.$ Then
$$
\mathcal F (\mathbf{f_{\alpha}}) \stackrel{S'(\mathbb{R}^n)}{=} C(\alpha, n) \cdot \mathbf{f}_{n - \alpha},
$$
where
$$
C(\alpha, n) = \frac{2^{\frac{n}{2} - \alpha} \cdot \Gamma(\frac{n -\alpha}{2})}{\Gamma(\frac{\alpha}{2})}.
$$
\end{proposition}
Let us note that in \cite[\S5.1, Lemma 1]{Stein} value of the constant differs from $C(\alpha, n)$ but this is simply due to different multiplicative constants in Fourier transform's definition.

Let $x \in \mathbb{R}^n$. Then we define a function $j_x \colon \mathbb{R}^n \longrightarrow \mathbb{C}$ by
$$
j_x(y) = e^{-i \: <x, \: y>} \; \; \forall \: y \in \mathbb{R}^n.
$$

Next lemma follows immediately from Proposition \ref{Fourier_transform} and equality
$$
(\mathcal{F}(\varphi))(x - t) = (\mathcal{F}^{-1}(j_x \cdot \varphi))(t) \; \; \; \forall \: t \in \mathbb{R}^n,
$$
which holds true for all $x \in \mathbb{R}^n$.

\begin{lemma}\label{convolution_Fourier_transform}
Let $0 < \alpha < n, \; \varphi \in S(\mathbb{R}^n)$. Then $\varphi \cdot f_{\alpha} \in L_1(\mathbb{R}^n)$ and the following equality holds
$$
\bigl(\mathcal{F}(\varphi \cdot f_{\alpha})\bigr)(x) = (2 \pi)^{-\frac n 2} \cdot C(\alpha, n) \cdot \bigl(\mathcal{F}(\varphi) \ast f_{n - \alpha}\bigr)(x) \; \; \forall \: x \in \mathbb{R}^n,
$$
where $\varphi \ast \psi$ is a convolution of functions $\varphi$ and $\psi.$
\end{lemma}

Our main technical tool in the proof of this section's main result is the following uniform estimate for the Fourier transform of $\psi_{(z)} \cdot f_{\alpha}$.

\begin{lemma}\label{estimate_Fourier_transform}
Let $\alpha \in (0; \: n)$ and $\psi \in D(\mathbb{R}^n)$. Then there exists a constant $M(\alpha, n, \psi) > 0$, such that for arbitrary $z \in \mathbb{R}^n$ we have
\begin{equation}\label{Fourier_transform_uniform_estimate}
| (\mathcal{F}(\psi_{(z)} \cdot f_\alpha))(x) | \leqslant M(\alpha, n, \psi) \cdot (1 + |x|^2)^{\frac{\alpha - n}{2}} \quad \forall \: x \in \mathbb{R}^n.
\end{equation}
\end{lemma}

Proof. Let us show that in order to prove \eqref{Fourier_transform_uniform_estimate} it is sufficient to establish estimates
\begin{equation}\label{F_T_eq_1}
|(\mathcal{F}(\psi_{(z)} \cdot f_\alpha))(x)| \leqslant M_1 \quad \forall \; x \in \mathbb{R}^n: |x| < 1
\end{equation}
and
\begin{equation}\label{F_T_eq_2}
|(\mathcal{F}(\psi_{(z)} \cdot f_\alpha))(x)| \leqslant M_2 \cdot |x|^{\alpha - n} \quad \forall \; x \in \mathbb{R}^n: |x| \geqslant 1
\end{equation}
for some constants $M_1 = M_1(\alpha, n, \psi) > 0$ and $M_2 = M_2(\alpha, n, \psi) > 0$, independent of $z$. Indeed, from \eqref{F_T_eq_1} and \eqref{F_T_eq_2} it follows that \eqref{Fourier_transform_uniform_estimate} is valid with a constant $M = 2^{\frac{n - \alpha}{2}}(M_1 + M_2)$, since for arbitrary $x \in B_1$ we obtain
$$
M_1 < \frac {2^{\frac{n - \alpha}{2}} \cdot (M_1 + M_2)}{(1 + |x|^2)^{\frac{n - \alpha}{2}}}
$$
and for arbitrary $x \in \mathbb{R}^n \setminus B_1$ we obtain
$$
\frac {M_2} {|x|^{n - \alpha}} \leqslant \frac{2^{\frac{n - \alpha}{2}} \cdot (M_1 + M_2)}{(1 + |x|^2)^{\frac{n - \alpha}{2}}} \: .
$$

Let us remark that we have a relation
$$
(\mathcal{F}(\psi_{(z)}))(x) = e^{-i \: < x, \: z >} \cdot (\mathcal{F}(\psi))(x) \; \; \forall \: x \in \mathbb{R}^n.
$$
Since from Lemma \ref{convolution_Fourier_transform} it follows that
$$
\mathcal{F}(\psi_{(z)} \cdot f_{\alpha}) = (2 \pi)^{-\frac n 2} \cdot C(\alpha, n) \cdot \mathcal{F}(\psi_{(z)}) * f_{n - \alpha},
$$
we obtain the following estimate:
$$
| (\mathcal{F}(\psi_{(z)} \cdot f_{\alpha}))(x) | \leqslant (2 \pi)^{-\frac n 2} \cdot C(\alpha, n) \cdot \int\limits_{\mathbb{R}^n} \frac{|(\mathcal{F}(\psi))(x - y)|}{|y|^{n - \alpha}} dy.
$$
Since for all $\psi \in D(\mathbb{R}^n) \subset S(\mathbb{R}^n)$ we have $\mathcal{F}(\psi) \in S(\mathbb{R}^n)$, convergence of the improper integral in the right--hand side of the estimate above follows from Proposition \ref{correctness_prop}.

As $\mathcal{F}(\psi) \in S(\mathbb{R}^n)$, it is obvious that
$$
\| \mathcal{F}(\psi) \|_{0, N} \stackrel{def}{=} \sup\limits_{x \in \mathbb{R}^n} \left((1 + |x|^2)^N \cdot |(\mathcal{F}(\psi))(x)|\right) < +\infty \; \; \forall \: N \in \mathbb{N},
$$
where $\| \cdot \|_{0, N}$ is one of the seminorms, which generate topology of Schwartz space $S(\mathbb{R}^n)$.

For arbitrary $N \in \mathbb{N}$ we denote constant $\| \mathcal{F}(\varphi) \|_{0, N}$ by $C_1(N, \varphi)$.

In a partial case when $N_1 = [\frac{n}{2}] + 1$ we arrive at inequality
$$
|(\mathcal{F}(\psi))(x - y)| \leqslant \: C_1(N_1, \psi) \cdot (1 + |x - y|^2)^{-N_1} \; \; \forall \: x, \: y \in \mathbb{R}^n.
$$
Hence,
\begin{equation}\label{F_T_eq_3}
| (\mathcal{F}(\psi_{(z)} \cdot f_{\alpha}))(x) | \leqslant (2 \pi)^{-\frac n 2} \cdot C(\alpha, n) \cdot C_1(N_1, \psi) \cdot \int\limits_{\mathbb{R}^n} \frac{dy}{|y|^{n - \alpha} \: (1 + |x - y|^2)^{N_1}} \: ,
\end{equation}
where convergence of improper integral in the right--hand side follows from the estimates $n - \alpha < n$ and $\: n - \alpha + 2 \cdot N_1 > 2 \cdot N_1 > n$.

In order to get an estimate of the integral in the right--hand side of \eqref{F_T_eq_3}, we treat cases $x \in \overline{B_1}$ and $x \notin \overline{B_1}$ separately.

Firstly, let us consider the case $x \in \overline{B_1}$, i.e. $|x| \leqslant 1$.

In this case we divide integral over $\mathbb{R}^n$ into the sum of integrals over $B_1$ and its complement and obtain the following estimate:
$$
\int\limits_{\mathbb{R}^n} \frac{dy}{|y|^{n - \alpha} (1 + |x - y|^2)^{N_1}} = \int\limits_{B_1} \frac{dy}{|y|^{n - \alpha} (1 + |x - y|^2)^{N_1}} + \int\limits_{\mathbb{R}^n \setminus B_1} \frac{dy}{|y|^{n - \alpha} (1 + |x - y|^2)^{N_1}} \leqslant
$$
$$
\leqslant \int\limits_{B_1} \frac{dy}{|y|^{n - \alpha}} + \int\limits_{\mathbb{R}^n \setminus B_1} \frac{dy}{|y|^{n - \alpha} (1 + |x - y|^2)^{N_1}}.
$$

Convergence of the integral $\int\limits_{B_1} \frac{dy}{|y|^{n - \alpha}}$ follows from the integrability criterion for $f_{\alpha}$ over unit ball $B_1$.

In order to estimate the integral $\int\limits_{\mathbb{R}^n \setminus B_1} \frac{dy}{|y|^{n - \alpha} (1 + |x - y|^2)^{N_1}}$ let us note that for arbitrary $y \notin B_1$ we have
$$
| x - y | \geqslant | y | - | x | \geqslant | y | - 1 \geqslant 0
$$
and thus
$$
| x - y |^2 \geqslant (| y | - 1)^2.
$$
Therefore, we obtain the estimate
$$
\int\limits_{\mathbb{R}^n \setminus B_1} \frac{dy}{|y|^{n - \alpha} \: (1 + |x - y|^2)^{N_1}} \leqslant \int\limits_{\mathbb{R}^n \setminus B_1} \frac{dy}{|y|^{n - \alpha} \: (1 + (|y| - 1)^2)^{N_1}} \; \; ,
$$
where the integral in the right--hand side of the last inequality doesn't depend on $x$ anymore and, as $n - \alpha + 2 \cdot N_1 > n$, its convergence readily follows from the integrability criterion for function $f_{\alpha}$ over $\mathbb{R}^n \setminus B_1$.

Thus, for any $x \in \overline{B_1}$ we obtain the estimate
$$
\int\limits_{\mathbb{R}^n} \frac{dy}{|y|^{n - \alpha} (1 + |x - y|^2)^{N_1}} \leqslant C_2(N_1, n, \alpha)
$$
for some constant $C_2(N_1, n, \alpha) > 0$.

Now let us treat the case $x \notin \overline{B_1}$. Here we also split the integral into two summands, yet using a different partition of $\mathbb{R}^n$:
$$
\int\limits_{\mathbb{R}^n} \frac{dy}{|y|^{n - \alpha} (1 + |x - y|^2)^{N_1}} = \int\limits_{B_{|x|/2}} \frac{dy}{|y|^{n - \alpha} (1 + |x - y|^2)^{N_1}}  \; + \int\limits_{\mathbb{R}^n \setminus B_{|x|/2}} \frac{dy}{|y|^{n - \alpha} (1 + |x - y|^2)^{N_1}}.
$$

Since from $|y| < \frac{|x|}{2}$ it follows that $|x - y| \geqslant |x| - |y| > \frac{|x|}{2}$, the first integral is estimated as follows:
$$
\int\limits_{B_{|x|/2}} \frac{dy}{|y|^{n - \alpha} (1 + |x - y|^2)^{N_1}} \leqslant \frac{1}{\left( 1 + \frac{|x|^2}{4} \right)^{N_1}} \int\limits_{B_{|x|/2}} \frac{dy}{|y|^{n - \alpha}} = \frac{C(n, \alpha) \cdot 4^{N_1}}{\left( 4 + |x|^2 \right)^{N_1}} \int\limits^{\frac{|x|}{2}}_0 \frac{r^{n - 1}}{r^{n - \alpha}} dr =
$$
$$
= \frac{C(n, \alpha) \cdot 2^{2 N_1}}{\alpha \left( 4 + |x|^2 \right)^{N_1}} \cdot \frac{|x|^{\alpha}}{2^{\alpha}} \leqslant \frac{C(n, \alpha) \cdot 2^{2 N_1 - \alpha}}{\alpha} \cdot \frac{1}{|x|^{2 N_1 - \alpha}} \leqslant D_1(N_1, n, \alpha) \cdot |x|^{\alpha - n},
$$
where $D_1(N_1, n, \alpha)$ is a positive constant and the last inequality is valid since $|x| > 1$ and $N_1 > \frac{n}{2}.$

For the second integral we have
$$
\int\limits_{\mathbb{R}^n \setminus B_{|x|/2}} \frac{dy}{|y|^{n - \alpha} (1 + |x - y|^2)^{N_1}} \leqslant \frac{2^{n - \alpha}}{|x|^{n - \alpha}} \int\limits_{\mathbb{R}^n \setminus B_{|x|/2}} \frac{dy}{(1 + |x - y|^2)^{N_1}} \leqslant
$$
$$
\leqslant 2^{n - \alpha} \cdot |x|^{\alpha - n} \int\limits_{\mathbb{R}^n} \frac{dy}{(1 + |x - y|^2)^{N_1}} = 2^{n - \alpha} \cdot |x|^{\alpha - n} \int\limits_{\mathbb{R}^n} \frac{du}{(1 + |u|^2)^{N_1}} = D_2(N_1, n, \alpha) \cdot |x|^{\alpha - n},
$$
where $D_2(N_1, n, \alpha)$ is a positive constant and the last inequality is valid because the integral under consideration converges for $2 \cdot N_1 > n.$

Thereby, for all $x \in \overline{B_1}$ we arrived at the estimate
$$
\int\limits_{\mathbb{R}^n} \frac{dy}{|y|^{n - \alpha} (1 + |x - y|^2)^{-N_1}} \leqslant C_2(N_1, n, \alpha),
$$
while for all $x \notin \overline{B_1}$ the estimate
$$
\int\limits_{\mathbb{R}^n} \frac{dy}{|y|^{n - \alpha} (1 + |x - y|^2)^{N_1}} \leqslant (D_1(N_1, n, \alpha) + D_2(N_1, n, \alpha)) \cdot |x|^{\alpha - n}
$$
was obtained.

From these estimates, alongside with the estimate \eqref{F_T_eq_3}, we get desired estimates \eqref{F_T_eq_1} and \eqref{F_T_eq_2}.

This concludes the proof of Lemma \ref{estimate_Fourier_transform}.

\bigskip

Now, using Lemma \ref{estimate_Fourier_transform}, we establish a sufficient condition for distribution $\mathbf{f_{\alpha}}$ to belong to the space $H^{-t}_{2, \: unif}(\mathbb{R}^n)$.

\medskip

\begin{proposition}\label{unif_sufficiency}
Let $t > - \frac{n}{2}.$ Then for arbitrary $\alpha \in \bigl(0; \: \min \bigl( n, t + \frac n 2 \bigr)\bigr)$
$$
\mathbf{f_{\alpha}} \in H^{-t}_{2, \: unif}(\mathbb{R}^n).
$$
\end{proposition}

Proof. Let $\eta$ be a function from $D(\mathbb{R}^n)$ such that
$$
0 \leqslant \eta(x) \leqslant 1 \;\; \forall \; x \in \mathbb{R}^n, \; \; \eta(x) = 1 \; \; \forall \; x \in \mathbb{R}^n: |x| \leqslant 1 \; \; \mbox{and} \; \; \eta(x) = 0 \; \; \forall \: x \in \mathbb{R}^n: |x| \geqslant 2.
$$

We need to prove that
$$
\eta_{(z)} \cdot \mathbf{f_{\alpha}} \in H^{-t}_2(\mathbb{R}^n) \; \; \; \forall \: z \in \mathbb{R}^n
$$
and obtain an upper estimate of $\sup\limits_{z \in \mathbb{R}^n} \| \eta_{(z)} \cdot \mathbf{f_{\alpha}} \|_{H^{-t}_2(\mathbb{R}^n)}$.
Since Fourier transform  is an isometric automorphism of the Banach space $(H^0_2(\mathbb{R}^n), \: \| \cdot \|_{H^0_2(\mathbb{R}^n)})$, this is equivalent to the fact that
$$
\mathcal{F}(J_{-t}(\eta_{(z)} \cdot \mathbf{f_{\alpha}})) \in H^0_2(\mathbb{R}^n) \; \; \forall \: z \in \mathbb{R}^n
$$
and there exists an upper estimate for
$$
\sup\limits_{z \in \mathbb{R}^n} \| \mathcal{F}(J_{-t}(\eta_{(z)} \cdot \mathbf{f_{\alpha}})) \|_{H^0_2(\mathbb{R}^n)} = \sup\limits_{z \in \mathbb{R}^n} \| \varphi_{-t} \cdot \mathcal{F}(\eta_{(z)} \cdot \mathbf{f}_{\alpha}) \|_{H^0_2(\mathbb{R}^n)}, 
$$
where
$$
\varphi_{-t}(x) = (1 + |x|^2)^{-\frac{t}{2}} \; \; \forall \: x \in \mathbb{R}^n.
$$


From $t > -\frac{n}{2}$ and the fact that $\alpha < n$ implies $\mathbf{f}_{\alpha} \in S'(\mathbb{R}^n)$, it can be deduced that $\varphi_{-t} \cdot \mathcal{F}(\eta_{(z)} \cdot \mathbf{f}_{\alpha})$ is well--defined as an element of $S'(\mathbb{R}^n)$. Since $\eta_{(z)} \cdot f_{\alpha} \in L_1(\mathbb{R}^n)$, we conclude that a distribution $\varphi_{-t} \cdot \mathcal{F}(\eta_{(z)} \cdot \mathbf{f}_{\alpha}) \in S'(\mathbb{R}^n)$ is a regular functional, generated by a function $\varphi_{-t} \cdot \mathcal{F}(\eta_{(z)} \cdot f_{\alpha}).$

Estimate $\alpha < n$ allows us to apply Lemma \ref{estimate_Fourier_transform}, which yields inequality
$$
(1 + |x|^2)^{-t} | (\mathcal{F}(\eta_{(z)} \cdot f_{\alpha}))(x)|^2 \leqslant \frac{(M(\alpha, n, \eta))^2}{(1 + |x|^2)^{t + n - \alpha}} \; \; \; \forall \: x \in \mathbb{R}^n.
$$
Hence, we deduce an uniform estimate
\begin{equation}\label{suff_estimate}
\int\limits_{\mathbb{R}^n} (1 + |x|^2)^{-t} | (\mathcal{F}(\eta_{(z)} \cdot f_{\alpha}))(x)|^2 \; dx \; \leqslant (M(\alpha, n, \eta))^2 \int\limits_{\mathbb{R}^n} \frac{dx}{(1 + |x|^2)^{t + n - \alpha}} \: ,
\end{equation}
where in the right--hand side integral converges as we have $\alpha < t + \frac{n}{2}$, which, in turn, implies $2 \cdot (t + n -\alpha) > n$.

Therefore, for arbitrary $z \in \mathbb{R}^n$ it was proved that
$$
\varphi_{-t} \cdot \mathcal{F}(\eta_{(z)} \cdot f_{\alpha}) \in L_2(\mathbb{R}^n)
$$
and, consequently,
$$
\eta_{(z)} \cdot \mathbf{f_{\alpha}} \in H^{-t}_2(\mathbb{R}^n).
$$
Since \eqref{suff_estimate} is valid, we thus proved that
$$
\sup\limits_{z \in \mathbb{R}^n} \| \eta_{(z)} \cdot \mathbf{f_{\alpha}} \|^2_{H^{-t}_2(\mathbb{R}^n)} \: = \: \sup\limits_{z \in \mathbb{R}^n} \int\limits_{\mathbb{R}^n} (1 + |x|^2)^{-t} | (\mathcal{F}(\eta_{(z)} \cdot f_{\alpha}))(x)|^2 \: dx < + \infty.
$$

But the latter estimate means that $\mathbf{f_{\alpha}} \in H^{-t}_{2, \: unif}(\mathbb{R}^n)$.

This concludes the proof of Proposition \ref{unif_sufficiency}.

\bigskip

\bigskip

\bigskip

{\Large 4. Optimal character of the embedding of uniformly localized Bessel potential space into multiplier space}

\bigskip
\medskip

In this section we shall establish necessary (Theorem \ref{multipliers_necessity}) and sufficient (Theorem \ref{multipliers_sufficiency}) conditions for the functional $\mathbf{f_{\alpha}}$ to be a multiplier from $H^s_2(\mathbb{R}^n)$ to $H^{-t}_2(\mathbb{R}^n)$. The necessary condition allows us to demonstrate sharp character of the index $\frac{n}{\max(s, t)}$ in the continuous embedding
$$
H^{-\min(s, t)}_{\frac{n}{\max(s, t)}, \: unif}(\mathbb{R}^n) \subset M[s, -t],
$$
which is the main result of our paper.

\begin{theorem}\label{multipliers_sufficiency}
Let $s, \: t \geqslant 0$ and $\max(s, t) < \frac{n}{2}$. Then
$$
\mathbf{f_{\alpha}} \in M[s, -t] \; \; \forall \: \alpha \in (0; \: s + t).
$$
\end{theorem}

Proof. First of all, since $M[s, -t] = M[t, -s]$ because of Remark \ref{multipliers_symmetry}, without loss of generality we may assume $0 \leqslant t \leqslant s < \frac{n}{2}$.

Let $s_1 = s + t - \frac{n}{2}$. Then
$$
s_1 > - \frac{n}{2} \; \; \; \mbox{and} \; \; \; \alpha < s + t = s_1 + n/2.
$$
Since estimates $\max(s, t) < \frac{n}{2}$ and $\alpha < s + t$ imply $\alpha < n$, by Lemma \ref{unif_sufficiency}, we have
$$
\mathbf{f}_{\alpha} \in H^{-s_1}_{2, \: unif}(\mathbb{R}^n).
$$

Since
$$
2 \leqslant \frac{n}{s} \; \; \; \; \mbox{и} \; \; -s_1 - \frac{n}{2} = -s - t = -t - \frac{n}{\frac{n}{s}} \: ,
$$
conditions of Remark \ref{unif_embedding} are met and, consequently, we get continuous embedding
$$
H^{-s_1}_{2, \: unif}(\mathbb{R}^n) \subset H^{-t}_{\frac{n}{s}, \: unif}(\mathbb{R}^n).
$$

Previously cited theorem from \cite{NZSh2} implies that continuous embedding
$$
H^{-t}_{\frac{n}{s}, \: unif}(\mathbb{R}^n) \subset M[s, -t]
$$
is valid and, consequently, we have
$$
\mathbf{f}_{\alpha} \in H^{-s_1}_{2, \: unif}(\mathbb{R}^n) \subset H^{-t}_{\frac{n}{s}, \: unif}(\mathbb{R}^n) \subset M[s, -t].
$$

This concludes the proof of Theorem \ref{multipliers_sufficiency}.

\begin{remark}\label{sufficiency_remark}
Analogous sufficient condition for the functional $f_{\alpha}$ to be a multiplier from $H^s_2(\mathbb{R}^n)$ to $H^{-t}_2(\mathbb{R}^n)$ was obtained in \cite[Lemma 6]{LRG} by other methods in a slightly less general situation when $s > t$.
\end{remark}

\begin{theorem}\label{multipliers_necessity}
Let $s, \: t \geqslant 0, \; \max(s, t) < \frac{n}{2}, \; 0 < \alpha < n$. Then $\mathbf{f_{\alpha}} \in M[s, -t]$ implies $\alpha \leqslant s + t.$
\end{theorem}

Proof. Let $f \in M[s, -t]$. Then, by definition of multiplier space, there exists a constant $C > 0$, such that
\begin{equation}\label{M_N_eq_1}
\mid \mathbf{f_{\alpha}}(g \cdot \overline{h}) \mid \; \leqslant C \: \| \mathbf{g} \|_{H^s_2(\mathbb{R}^n)} \| \mathbf{h} \|_{H^t_2(\mathbb{R}^n)} \qquad \forall \; g, \: h \in D(\mathbb{R}^n).
\end{equation}

Since for arbitrary $\gamma \in \mathbb{R}$ and $p > 1 \;$ $D(\mathbb{R}^n)$ is dense in the space $(S(\mathbb{R}^n), \| \cdot \|_{\mathbb{H}^{\gamma}_p(\mathbb{R}^n)})$ (where $\| f \|_{\mathbb{H}^{\gamma}_p(\mathbb{R}^n)} \stackrel{def}{=} \| \mathbf{f} \|_{H^{\gamma}_p(\mathbb{R}^n)}$), Schwartz space $S(\mathbb{R}^n)$ is a topological algebra with respect to pointwise multiplication and $\mathbf{f}_{\alpha}$ belongs to $S'(\mathbb{R}^n)$ for any $\alpha \in (0; \: n)$ estimate \eqref{M_N_eq_1} can be extended onto the whole space $S(\mathbb{R}^n) \colon$
\begin{equation}\label{M_N_eq_2}
\mid \mathbf{f_{\alpha}}(g \cdot \overline{h}) \mid \; \leqslant C \: \| \mathbf{g} \|_{H^s_2(\mathbb{R}^n)} \| \mathbf{h} \|_{H^t_2(\mathbb{R}^n)} \qquad \forall \; g, \: h \in S(\mathbb{R}^n).
\end{equation}

Since for arbitrary $\gamma \in \mathbb{R}$ we have
$$
\| \mathbf{f} \|_{H^{\gamma}_2(\mathbb{R}^n)} = \| J_{\gamma}(f) \|_{L_2(\mathbb{R}^n)} = \| \mathcal{F}(J_{\gamma}(f)) \|_{L_2(\mathbb{R}^n)} = \| \varphi_{\gamma} \cdot \mathcal{F}(f) \|_{L_2(\mathbb{R}^n)} \; \; \forall \; f \in S(\mathbb{R}^n),
$$
and $g \cdot f_{\alpha} \in L_1(\mathbb{R}^n)$ implies that
$$
\mathbf{f_{\alpha}}(g \cdot \overline{h}) = (\mathcal{F}(g \cdot \mathbf{f_{\alpha}})) (\mathcal{F}^{-1}(\overline{h})) = (\mathcal{F}(g \cdot \mathbf{f_{\alpha}}))(\overline{\mathcal{F}(h)}),
$$
inequality \eqref{M_N_eq_2} can be rewritten as
\begin{equation}\label{M_N_eq_3}
|(\mathcal{F}(g \cdot \mathbf{f_{\alpha}}))(\overline{\mathcal{F}(h)})| \leqslant \; C \: \| \varphi_s \cdot \mathcal{F}(g) \|_{L_2(\mathbb{R}^n)} \: \| \varphi_t \cdot \mathcal{F}(h) \|_{L_2(\mathbb{R}^n)} \; \; \; \forall \: g, \: h \in S(\mathbb{R}^n).
\end{equation}

Let us recall that a convolution of a function from the Schwartz space $S(\mathbb{R}^n)$ and a distribution from the dual Schwartz space $S'(\mathbb{R}^n)$ is defined as follows:
$$
(\psi \ast u)(f) = u(\psi_{-} \ast f) \; \; \; \forall \: u \in S'(\mathbb{R}^n), \; \forall \: \psi, \: f \in S(\mathbb{R}^n),
$$
where $\psi_{-} \colon \mathbb{R}^n \to \mathbb{C}, \; \; x \mapsto \psi(-x).$

Since $g \in S(\mathbb{R}^n), \; \mathbf{f_{\alpha}} \in S'(\mathbb{R}^n),$ and, as it is well-known (see, e.g., \cite[Theorem 5.18]{GrBook}) that
$$
\mathcal{F}(\psi \cdot u) = (2 \cdot \pi)^{\frac{n}{2}} \bigl(\mathcal{F}(\psi)\bigr) \ast \bigl(\mathcal{F}(u)\bigr) \; \; \; \forall \: \psi \in S(\mathbb{R}^n), \; \forall \; u \in S'(\mathbb{R}^n),
$$
we arrive at the following equality
$$
(\mathcal{F}(g \cdot \mathbf{f_{\alpha}}))(\overline{\mathcal{F}(h)}) = (2 \cdot \pi)^{\frac{n}{2}} \cdot C(\alpha, n) \cdot \bigl( \mathcal{F}(g) \ast \mathbf{f}_{n - \alpha} \bigr)(\overline{\mathcal{F}(h)}) =
$$
$$
= (2 \cdot \pi)^{\frac{n}{2}} \cdot C(\alpha, n) \cdot \mathbf{f}_{n - \alpha} \bigl((\mathcal{F}(g))_{-} \ast \overline{\mathcal{F}(h)} \bigr).
$$

Employing this equality, it follows immediately that \eqref{M_N_eq_3} is equivalent to inequality
\begin{equation}\label{M_N_eq_4}
\bigl| \mathbf{f}_{n - \alpha} \bigl((\mathcal{F}(g))_{-} \ast \overline{\mathcal{F}(h)} \bigr) \bigr| \; \leqslant \; K(\alpha, n) \: \| \varphi_s \cdot \mathcal{F}(g) \|_{L_2(\mathbb{R}^n)} \: \| \varphi_t \cdot \mathcal{F}(h) \|_{L_2(\mathbb{R}^n)} \; \; \; \forall \: g, \: h \in S(\mathbb{R}^n),
\end{equation}
where
$$
K(\alpha, n) = (2 \cdot \pi)^{-\frac{n}{2}} \cdot \frac{C}{C(\alpha, n)} \: .
$$

Now for arbitrary function $\psi \in S(\mathbb{R}^n)$ and arbitrary $r \geqslant 0$ let us define function $\widetilde{\psi_r} \colon \mathbb{R}^n \to \mathbb{R}$ by letting
$$
\widetilde{\psi_r} = \varphi_r \cdot (\mathcal{F}(\psi)).
$$
Since both Fourier transform and multiplication operator $A_{\varphi_r}$ by a function $\varphi_r, \: \gamma \in \mathbb{R}$, are linear homeomorphisms of $S(\mathbb{R}^n)$ onto itself, we see that \eqref{M_N_eq_4} holds true if and only if the following inequality is valid
$$
\Bigl| \: \int\limits_{\mathbb{R}^n} \biggl( \int\limits_{\mathbb{R}^n} \varphi_{-s}(y - x) \cdot g_0(y - x) \cdot \varphi_{-t}(y) \cdot \overline{h_0(y)} dy \biggr) \cdot f_{n - \alpha}(x) dx \Bigr| \leqslant \; K(\alpha, n) \: \| g_0 \|_{L_2(\mathbb{R}^n)} \| h_0 \|_{L_2(\mathbb{R}^n)}
$$
for arbitrary functions $g_0, \: h_0 \in S(\mathbb{R}^n)$.

By Fubini's theorem, we can change the order of integration and, taking into account the relation
$$
(\varphi_{-s} \cdot g_0) \ast f_{n - \alpha} = f_{n - \alpha} \ast (\varphi_{-s} \cdot g_0),
$$
we rewrite the latter inequality as
\begin{equation}\label{M_N_eq_5}
\Bigl| \: \int\limits_{\mathbb{R}^n} \biggl( \; \int\limits_{\mathbb{R}^n} f_{n - \alpha}(y - x) \cdot \varphi_{-s}(x) \cdot g_0(x) dx \biggr) \; \varphi_{-t}(y) \cdot \overline{h_0(y)} \: dy \Bigr| \; \leqslant \; K(\alpha, n) \: \| g_0 \|_{L_2(\mathbb{R}^n)} \| h_0 \|_{L_2(\mathbb{R}^n)}.
\end{equation}

For arbitrary $m \in \mathbb{N}$ we define a function $\eta_m$ by letting
$$
\eta_m(x) = \eta\Bigl(\frac{x}{m}\Bigr) \; \; \forall \: x \in \mathbb{R}^n,
$$
where $\eta \in D(\mathbb{R}^n) \:$ --- is a real--valued function satisfying conditions
$$
0 \leqslant \eta(x) \leqslant 1 \; \; \forall \: x \in \mathbb{R}^n, \; \; \; \eta(x) = 1 \; \; \forall \: x \in \mathbb{R}^n \colon |x| \leqslant 1, \; \; \; \eta(x) = 0 \; \; \forall \: x \in \mathbb{R}^n \colon |x| \geqslant 2.
$$
Taking $g_0 = h_0 = \eta_m$ in the inequality \eqref{M_N_eq_5}, we obtain estimate
$$
\Bigl| \: \int\limits_{\mathbb{R}^n} \biggl( \; \int\limits_{\mathbb{R}^n} f_{n - \alpha}(y - x) \cdot \varphi_{-s}(x) \cdot \eta \left( \frac{x}{m} \right) dz \biggr) \: \varphi_{-t}(y) \cdot \eta \left( \frac{y}{m} \right) dy \Bigr| \; \leqslant \; K(\alpha, n) \: \| \eta_m \|_{L_2(\mathbb{R}^n)}^2,
$$
for arbitrary $m \in \mathbb{N}$.

Since all functions in the left--hand side of this estimate are non--negative, we get
$$
\Bigl| \: \int\limits_{\mathbb{R}^n} \biggl( \; \int\limits_{\mathbb{R}^n} f_{n - \alpha}(y - x) \cdot \varphi_{-s}(x) \cdot \eta \left( \frac{x}{m} \right) dx \biggr) \: \varphi_{-t}(y) \cdot \eta \left( \frac{y}{m} \right) dy \Bigr| \geqslant
$$
$$
\geqslant \int\limits_{B_m(0)} \biggl( \; \int\limits_{B_m(0)} f_{n - \alpha}(y - x) \cdot \varphi_{-s}(x) \cdot \eta \left( \frac{x}{m} \right) dx \biggr) \: \varphi_{-t}(y) \cdot \eta \left( \frac{y}{m} \right) dy =
$$
$$
= \int\limits_{B_m(0)} \biggl( \; \int\limits_{B_m(0)} |y - x|^{-(n - \alpha)} (1 + |x|^2)^{-\frac{s}{2}} dx \biggr) \: (1 + |y|^{2})^{-\frac{t}{2}} dy.
$$
Also we have
$$
\| \eta_m \|_{L_2(\mathbb{R}^n)}^2 = \int\limits_{\mathbb{R}^n} (\eta(x/m))^2 dx = m^n \cdot \int\limits_{\mathbb{R}^n} (\eta(z))^2 dz = C_1 \cdot m^n,
$$
where $C_1 = \| \eta \|^2_{L_2(\mathbb{R}^n)}$.

Applying inequality
$$
\frac{1}{|y - x|^{n - \alpha} (1 + |x|^2)^{\frac{s}{2}}} \geqslant \frac{1}{[2(1 + |x|^{2} + |y|^{2})]^{\frac{n - \alpha}{2}} (1 + |x|^{2})^{\frac{s}{2}}} \geqslant \frac{1}{2^{\frac{n - \alpha}{2}}(1 + |x|^{2} + |y|^{2})^{\frac{n - \alpha + s}{2}}} \; \; ,
$$
and changing from iterated integral to double integral we get
$$
\int\limits_{B_m(0) \times B_m(0)} \frac{(1 + |y|^{2})^{-\frac{t}{2}}}{(1 + |x|^{2} + |y|^{2})^{\frac{n - \alpha + s}{2}}} \; d(x \times y) \leqslant K(\alpha, n) \cdot C_1 \cdot 2^{\frac{n - \alpha}{2}} \cdot m^n.
$$
From this estimate we immediately obtain
\begin{equation}\label{M_N_eq_6}
\int\limits_{B_m(0) \times B_m(0)} \frac{1}{(1 + |x|^{2} + |y|^{2})^{\frac{n - \alpha + s + t}{2}}} \; d(x \times y) \leqslant K_1 \cdot m^n
\end{equation}
for the constant
$$
K_1 = K_1(\alpha, n, \eta) = K(\alpha, n) \cdot C_1 \cdot  2^{\frac{n - \alpha}{2}}.
$$

Since the ball $U_m(0) = \{(x, y) \in \mathbb{R}^n \times \mathbb{R}^n \: | \: |x|^{2} + |y|^{2} \leqslant m^{2} \}$ is a subset of $B_m(0) \times B_m(0)$, and integrand is a non--negative function, then, changing domain of integration from $B_m(0) \times B_m(0)$ to $U_m(0)$, we get
$$
K_1 \cdot m^n \; \geqslant \int\limits_{U_m(0)} \frac{1}{(1 + |x|^{2} + |y|^{2})^{\frac{n - \alpha + s + t}{2}}} \: d(x \times y) = C(2 \: n) \cdot \int\limits_{[0, \: m]} \frac{r^{2 n - 1} dr}{(1 + r^{2})^{\frac{n - \alpha + s + t}{2}}} \geqslant
$$
$$
\geqslant \: C(2 \: n) \cdot \int\limits_{[1, \: m]} \frac{r^{2 n - 1} dr}{(1 + r^{2})^{\frac{n - \alpha + s + t}{2}}} \geqslant C(2 \: n) \cdot \int\limits_{[1, \: m]} \frac{r^{2 n - 1} dr}{(2 \: r^{2})^{\frac{n - \alpha + s + t}{2}}} = C_2(n) \cdot \int\limits_{[1, m]} r^{n + \alpha - s - t - 1} dr,
$$
where
$$
C(2 n) = \frac{2 \cdot \pi^{n}}{\Gamma(n)} \: , \; \; \mbox{a} \; \; C_2(n) = C(2 n) \cdot 2^{\frac{\alpha - n - s - t}{2}}.
$$

Since $\max (s, t) < \frac{n}{2}$ implies $n + \alpha - s - t > 0$, finally we obtain
$$
m^{n + \alpha - s - t} < \frac{K_1 \cdot (n + \alpha - s - t)}{C_1(n)} \cdot m^n + 1.
$$
As this inequality is valid for arbitrary $m \in \mathbb{N}$, then we get estimate $n + \alpha - s - t \leqslant n$, from which our desired estimate $\alpha \leqslant s + t$ immediately follows.

This concludes the proof of Theorem \ref{multipliers_necessity}.

The principal result of the paper is the following theorem which demonstrates that the constant $p_1 = \frac{n}{\max(s, \; t)}$ is sharp in the sense that distracting arbitrarily small $\varepsilon > 0$ we no longer have embedding $H^{-\min(s, \: t)}_{p_1 - \varepsilon, \: unif} \subset M[s, \: -t]$.

\begin{theorem}\label{exact_embedding}
Let $s, \: t \geqslant 0$ and $0 < \max (s, t) < \frac{n}{2}.$ Then for arbitrary  $\varepsilon$, satisfying conditions
$$
0 < \varepsilon < \frac{n}{\max (s, t)} - 2,
$$
there exists a positive number
$$
\delta(\varepsilon) \in \Bigl(0; \frac{n}{2} - \max (s, t)\Bigr),
$$
such that for $\alpha = s + t + \delta(\varepsilon)$
$$
\mathbf{f_{\alpha}} \in H^{- \min (s, t)}_{\frac{n}{\max (s, t)} - \varepsilon, \: unif}(\mathbb{R}^n) \setminus \: M[s, -t].
$$
\end{theorem}

Proof. Proposition \ref{multipliers_symmetry} and symmetricity of Theorem \ref{exact_embedding} conditions allows us to assume, without loss of generality, that $0 \leqslant t \leqslant s < \frac{n}{2}$ and $s > 0.$

Fix an arbitrary number $\varepsilon$ such that
$$
0 < \varepsilon < \frac{n}{s} - 2 = \frac{n - 2 \: s}{s}
$$
and define
$$
\delta(\varepsilon) = \frac{s^2}{2 (\frac{n}{\varepsilon} - s)} \: .
$$
Since $0< \varepsilon < n/s - 2$, we have $n/\varepsilon - s > 0$ and, consequently, $\delta(\varepsilon) > 0$.

Let us check the condition
$$
\delta(\varepsilon) < \frac{n}{2} - s.
$$
Indeed,
$$
\frac{s^2}{2 (\frac{n}{\varepsilon} - s)} < \frac{n}{2} - s \; \Longleftrightarrow \; \frac{n}{\varepsilon} - s > \frac{s^2}{n - 2 s} \; \Longleftrightarrow \; \varepsilon < \frac{n}{n - s} \cdot \frac{n - 2 \: s}{s},
$$
where the last of these inequalities trivially follows from
$$
\varepsilon < \frac{n - 2 \: s}{s} \: .
$$

Then for $\alpha \stackrel{def}{=} s + t + \delta$ we have
$$
\alpha < s + t + \frac{n}{2} - s = t + \frac{n}{2} < n,
$$
hence, functional $\mathbf{f_{\alpha}}$ is well--defined as an element of $S'(\mathbb{R}^n).$

Letting
$$
t_1 \stackrel{def}{=} s + t + 2 \: \delta - \frac{n}{2} \: ,
$$
we obtain
$t_1 > -\frac{n}{2}$.

Since
$$
\alpha = s + t + \delta < s + t + 2 \: \delta = t_1 + \frac{n}{2} \: ,
$$
and  $\alpha < n$, then, by Theorem \ref{unif_sufficiency},
$$
\mathbf{f_{\alpha}} \in H^{-t_1}_{2, \: unif}(\mathbb{R}^n).
$$

From the chain of equalities
$$
t_1 + \frac{n}{2} = s + t + 2 \: \delta = t + s + \frac{s^2}{\frac{n}{\varepsilon} - s} = t + \frac{n \cdot s}{n - s \cdot \varepsilon}
$$
it follows that
$$
- t_1 - \frac{n}{2} = -t - \frac{n}{\frac{n}{s} - \varepsilon} \: .
$$
Since we also have $2 < n/s - \varepsilon$ we can apply the result of Remark \ref{unif_embedding}, which yelds a continuous embedding
$$
H^{-t_1}_{2, \: unif}(\mathbb{R}^n) \subset H^{-t}_{\frac{n}{s} - \varepsilon, \: unif}(\mathbb{R}^n)
$$

Therefore, we obtain
$$
\mathbf{f_{\alpha}} \in H^{-t_1}_{2, \: unif}(\mathbb{R}^n) \subset H^{-t}_{\frac{n}{s} - \varepsilon, \: unif}(\mathbb{R}^n).
$$
On the other hand, since $\alpha = s + t + \delta > s + t,$ then Theorem \ref{multipliers_necessity} implies
$$
\mathbf{f_{\alpha}} \notin M[s, -t].
$$

This concludes the proof of Theorem \ref{exact_embedding}.

\begin{remark}\label{unif_description_impossibility}
From the result of Theorem \ref{exact_embedding} it follows that in the case of $\max(s, t) < \frac{n}{2}$ constructive description of multiplier space $M[s, \: -t]$ in terms of the scale $H^{\gamma}_{p, \: unif}(\mathbb{R}^n)$ can not be established. Therefore, in this case we can not obtain an analogue of the result from \cite{BSh}, which states the coincidence of the spaces $M[s, -t]$ and $H^{-\min(s, t)}_{2, \: unif}(\mathbb{R}^n)$ and equivalence of their norms whether $\max(s, t) > \frac{n}{2}$.
\end{remark}

\bigskip

Addresses:

A.\,A.~Belyaev,
Lomonosov Moscow State University,
Department of Mechanics and Mathematics;
{email: alexei.a.belyaev@gmail.com}

\medskip
A.\,A.\, Shkalikov,
Lomonosov Moscow State University,
Department of Mechanics and Mathematics;
{email: shkalikov@mi.ras.ru}

\end{document}